\documentclass[pre]{revtex4}

\usepackage{amsmath}
\usepackage{amsfonts}
\usepackage{hyperref}
\usepackage{graphicx}

\usepackage[T1]{fontenc}

\renewcommand{\Re}{\mathrm{Re}}%
\renewcommand{\i}{\mathrm{i}}%

\begin{document}

\title{Factorial Series Representation of  Stieltjes Series Converging Factors}

\author{Riccardo Borghi}
\affiliation{Dipartimento di Ingegneria Civile, Informatica e delle Tecnologie Aeronautiche, {Universit\`a}  degli Studi ``Roma Tre'', {00146 Rome, Italy; } 
}

\begin{abstract}
The practical usefulness of Levin-type nonlinear sequence transformations as numerical tools
for the summation of divergent series or for the convergence acceleration of slowly 
converging series, is nowadays beyond dispute. Weniger's transformation, in particular, 
is able to accomplish spectacular results when used to overcome resummation problems, 
often outperforming better known resummation techniques, the most known being Pad\'e approximants. 
However, our understanding of its theoretical features is still far from being satisfactory and particularly bad as far as the decoding of factorially divergent series is concerned. 

Stieltjes series represent a class of power series of fundamental interest in mathematical physics.
In the present paper, it is shown how the Stieltjes series converging factor of any order is expressible
as an inverse factorial series, whose terms can be analytically retrieved through a simple recursive algorithm.  
A few examples of applications of our algorithm are presented, in order to show its effectiveness and implementation ease.
We believe the results presented here could constitute an important, preliminary step for the development of a general 
convergence theory of Weniger's transformation on Stieltjes series. 
A rather ambitious project, but worthy of being pursued in the future.
\end{abstract}

\maketitle

\section{Prelude}
\label{Sec:Prelude}

{\em
The present paper is a fitting tribute to the memory of Ernst Joachim Weniger, mentor and friend, who passed away on August 10th, 2022, sadly too soon. I had the good fortune to meet EJW in 2009 at a conference organized at the {\em Centre International de Rencontres Mathématiques} in Marseille, France. Since then, I have had the pleasure of working with him for over a decade. In 2015, we published our first (and sole) joint paper, in which a convergence theory for the resummation of the Euler series via Levin-type nonlinear sequence transformations, was presented. In particular, we proved that {\em his}
transformation, named Weniger transformation, were not only able to resum the Euler series~but also that, in accomplishing such task, it turned out to be ``exponentially faster''  than Pad\'e approximants. 
Euler series is an example of Stieltjes series, for which a well established general convergence theory based on Pad\'e approximants is well known to be available. After our 2015 joint paper, EJW and I thought that trying to conceive a convergence theory also for Levin-type transformations, when applied to Stieltjes series, could have been an  ambitious scientific project, but worthy of being pursued.

What is contained in the present work should have been my second joint work with EJW. In writing it, I tried to respect to the best of my ability all that EJW had and would suggest to me, also thanks to the help of still unpublished material and several discussions that took place during his visits to my former Engineering Department at the University "Roma Tre.''
}

\section{Introduction}
\label{Sec:Intro}

Power series are among the most important mathematical objects, not
only in classical analysis, but also in the mathematical treatment of several
scientific and engineering problems. 
Most power series have a \emph{finite} radius of convergence. If the argument of a power series is
close to the expansion point, convergence is usually (very) good, but once the argument approaches 
the boundary of the circle of convergence, convergence can become arbitrarily bad.  Outside the circle 
of convergence, power series diverge \cite{Hardy/1949}. 

Especially in theoretical physics, wildly divergent power series occur abundantly which cannot be replaced easily by alternative expressions. 
Classic examples of useful divergent series are the Stirling series for the logarithm of the gamma function ~\cite*[Eq.\
(5.11.1)]{Olver/Lozier/Boisvert/Clark/2010}, 
as well as 
the Poincar\'e seminal work on asymptotic series~\cite{Poincare/1886}, inspired by his earlier work in astronomy.
In quantum physics, divergent series have become indispensable. Dyson~\cite{Dyson/1952} argued that perturbation expansions
in quantum electrodynamics must diverge factorially. In their pioneering work on anharmonic oscillators, 
Bender and Wu ~\cite{Bender/Wu/1969,Bender/Wu/1971,Bender/Wu/1973} showed that factorially divergent
perturbation expansions occur also in nonrelativistic quantum mechanics. Later, it was found that factorially divergent perturbation
expansions are actually the rule in quantum physics rather than the exception (see for example the articles by ~\cite[Table 1]{Fischer/1997}
and ~\cite{Suslov/2005}, or the articles reprinted in the book by
~\cite{LeGuillou/Zinn-Justin/1990}).
In view of the conceptual and technical problems caused by divergent series both in applied mathematics and theoretical physics, it should not
be surprising that there has been an extensive research on summation techniques. Recent reviews including numerous references can be found in~\cite{Weniger/2008,Brezinski/RedivoZaglia/Weniger/2010a,Borghi/Weniger/2015,Borghi/2016}.

The most important summation techniques are Borel summation
\citep{Borel/1899}, which replaces a divergent perturbation expansion by
a Laplace-type integral, and Pad\'{e} approximants \citep{Pade/1892},
which transform the partial sums of a (formal) power series to rational
functions. There is an endless literature on these summation
techniques. Borel summation is discussed, for example, in \citet{Costin/2009}, \citet{Shawyer/Watson/1994},
\citet{Sternin/Shatalov/1996}, while the probably most complete treatment of
Pad\'{e} approximants can still be found in the monograph by
\citet{Baker/Graves-Morris/1996}.
Both Borel summation and Pad\'{e} approximants have -- like all other numerical techniques -- certain shortcomings and
limitations. In recent years, a considerable deal of work 
has been invested on the so-called sequence transformations, which turned out to be very useful numerical tools for
the summation of divergent series. 
An extended recent literature, in particular on \emph{nonlinear} and \emph{nonregular} sequence transformations, is available 
\cite{Brezinski/1977,Brezinski/1978,Brezinski/RedivoZaglia/1991a,Sidi/2003,Wimp/1981,Weniger/1989,Homeier/1996,Brezinski/RedivoZaglia/Weniger/2010a,Weniger/2010a,Aksenov/Savageau/Jentschura/Becher/Soff/Mohr/2003,Bornemann/Laurie/Wagon/Waldvogel/2004,Gil/Segura/Temme/2007,Caliceti/Meyer-Hermann/Ribeca/Surzhykov/Jentschura/2007,Temme/2007,Gil/Segura/Temme/2011,Trefethen/2013}.

Levin-type transformations \cite{Levin/1973,Weniger/1989,Weniger/2004}, in particular,
constitute  an  effective and powerful class of nonlinear sequence transformations, which 
differ substantially from other, better known sequence transformations, as for example the celebrated Wynn's
epsilon algorithm \citep{Wynn/1956a}.
Such a lot of successful applications of Levin-type transformations have been achieved and described in the literature (see for instance the references in ~\cite[Appendix B]{Weniger/2008},~\cite{Borghi/Weniger/2015},~\cite{Chang/He/Hu/Sun/Weniger/2020}, or in the articles 
in~\cite{Brezinski/RedivoZaglia/Weniger/2010a}), that they are also treated in NIST Handbook \cite[Chapter 3.9(v)
Levin's and Weniger's Transformations]{Olver/Lozier/Boisvert/Clark/2010}.
Now, from the perspective of who is interested only in employing Levin-type
transformations as computational tools, the situation looks quite good. 
The practical usefulness of sequence transformations as computational
tools for the acceleration of convergence as well as the summation of divergent
series is now established beyond doubt.
However, from a purely theoretical perspective, our understanding of the convergence properties
of Levin-type sequence transformations is still far from being satisfactory.

In the case of Pad\'{e}, the situation is much better. For example, if the input data
are the partial sums of a Stieltjes series, it can be
proved rigorously that certain subsequences of the Pad\'{e} table
converge to the value of the corresponding Stieltjes function.
More precisely, if a factorially divergent power series is a Stieltjes series
and if it also satisfies the so-called Carleman condition, 
such a series would be Pad\'{e} summable.
In ~\cite{Borghi/Weniger/2015}, 
a \emph{rigorous} convergence study of the summation of the
Euler series, namely $\mathcal{E} (z) \sim \sum_{n=0}^{\infty} (-1)^{n} n! z^{n}$, which is the paradigm of the factorially divergent power series
(for a compact description, see~\cite[Sec. 2]{Borghi/Weniger/2015}), was carried out by comparing the performances of Pad\'{e} approximants with those of a certain class of nonlinear sequence transformations.

Readers are encouraged to go through Section~4 of~\cite{Borghi/Weniger/2015}, where
an extensive  summary of the most relevant properties of Levin-type transformations can be found. 
Levin-type transformations \emph{always} admit explicit analytical expressions,
which might be used as a natural starting point for building up a convergence theory. 
However, there is a main conceptual drawback against such an ambitious project: 
the technical problems occurring in any convergence studies, substantially differ among the various Levin-type transformations. 
For example, in~\cite{Borghi/Weniger/2015} the 
convergence analysis of the summation of the Euler series turns out to be, in the case of the  Weniger transformation, much
easier than in the case of other Levin-type transformations.
The principal reason stems in the role played by the so-called {\em converging factor}~\cite{Airey/1937},\footnote{Converging factors must not be mixed up with
 \emph{convergence factors} for summable series, which are for instance discussed in an article and a book by~\cite{Moore/1907,Moore/2008}} 
whose theory has extensively been discussed in articles~\cite{Dingle/1958a,Dingle/1958b,Dingle/1958c,Dingle/1959a,Dingle/1959b,Dingle/1959c} 
and a book~\cite[Sections 21 -~26]{Dingle/1973} by \citeauthor{Dingle/1973}. 

Euler series  is a Stieltjes series, and its partial sums are Pad\'{e} summable to the so-called Euler integral.
In the present paper it is proved that, when the series under investigation is of Stieltjes type, an inverse factorial 
representation of its converging factor can always be analytically retrieved via an algorithm based of the following points:
\begin{itemize}

\item[(i)]
the converging factor of a typical Stieltjes series must satisfy a simple first-order finite difference equation;

\item[(ii)] 
there is a tight connection between inverse factorial series and recurrence relationships, 
the former being the natural mathematical tools for representing the solution of linear finite difference equations, 
similarly as solutions of differential equations are naturally expressed in terms of power series.  
A detailed discussion of the most important features of factorial series plus additional references can be found, for instance, in~\cite{Weniger/2010b} and a more condensed one in~\cite[Appendix B]{Borghi/Weniger/2015}. 
\end{itemize}

The paper is structured as follows: in Sec.~\ref{Sec:StieltjesSeriesConvFact}, after  the main definitions and properties
of Stieltjes series and Stieltjes functions have been briefly resumed, the recurrence rule for the converging factors of
Stieljes series is derived and formally solved through inverse factorial series expansions. 
In Sec.~\ref{SubSec:FactorialSeriesWenigerTransformation}, a, somewhat didactical derivation of Weniger's transformation is proposed, in order
to enlightened its connection with the factorial series representation of converging factors.
Section~\ref{SubSec:FactorialSeriesExpansion} represents the core of the paper: the algorithm to generate the factorial series representation of a typical
Stieltjes series converging factor is proposed and its derivation detailed.
Section~\ref{Sec:Examples} is devoted to illustrate several examples of application of the results established in 
Sec.~\ref{SubSec:FactorialSeriesExpansion}. Finally, a few conclusive considerations and plans for future works are outlined in Sec.~\ref{Sec:Conclusions}.

In order to improve paper's readability, the most technical parts of the paper have been confined within several appendices.


\section{Converging Factors of Stieltjes Series}
\label{Sec:StieltjesSeriesConvFact}
%



Consider an increasing real-valued  function $\mu(t)$ defined for $t\in[0,\infty]$, with~infinitely many points of increase. 
The measure $\mathrm{d}\mu$ is then positive on $[0,\infty]$. Assume all positive 	quantities 
\begin{equation}
\label{Eq:Stieltjes.0.1}
\displaystyle
\mu_m=\int_0^\infty\,t^m\,\mathrm{d}\mu\,,\qquad m \ge 0\,,
\end{equation}
which will be called moments, to be finite. Then, the~formal power series
\begin{equation}
\label{Eq:Stieltjes.0.2}
\displaystyle
\sum_{m=0}^\infty\,\dfrac{(-)^m}{z^{m+1}}\,\mu_m\,,
\end{equation}
is called a \emph{Stieltjes series}. Such series turns out to be asymptotic, for~ $z\to \infty$,  
to the function $F(z)$  defined by 
\begin{equation}
\label{Eq:Stieltjes.0.3}
\displaystyle
F(z)\,=\,\int_0^\infty\,\frac{\mathrm{d}\mu}{z+t}\,,\qquad |\arg(z)|<\pi\,,
\end{equation}
%
which is called \emph{Stieltjes function}. In other words,
\begin{equation}
\label{Eq:Stieltjes.3}
\displaystyle
F(z) \sim \sum_{m=0}^\infty\,\frac {(-)^m}{z^{m+1}}\,\mu_m\,,\qquad  z \to \infty\,.
\end{equation}
Whether such series converges or diverges depends on the behavior of the
moment sequence  $\{\mu_m\}^\infty_{m=0}$ as $m\to\infty$. 
Any Stieltjes function can always be expressed as the $n$th-order partial sum of the associated Stieltjes series (\ref{Eq:Stieltjes.3}),
plus a truncation error term which has itself the form of a Stieltjes integral (see for example~\cite[Theorem 13-1]{Weniger/1989}), i.e.,
\begin{equation}
  \label{Def_StieFunParSum}
  F (z) \, = \, \sum_{m=0}^{n} \, \dfrac{(-1)^{m}}{z^{m+1}} \, \mu_{m}  
  \, + \, \left(-\dfrac 1z\right)^{n+2} \, \int_{0}^{\infty} \, \frac{t^{n+1}\mathrm{d} \mu} {z+t} \, , \qquad \vert \arg (z) \vert < \pi \, .  
\end{equation}
Such a truncation error term can always be recast as the product of the first term neglected,
namely $(-1/z)^{n+2}\,\mu_{n+1}$, and a {\em converging factor}~\cite{Airey/1937,Dingle/1973}, 
say $\varphi_{n+1}(z)$, which is defined as follows:
\begin{equation}
\label{Eq:Stieltjes.7}
\displaystyle
\varphi_{m}(z) \,=\,\frac {1}{\mu_{m}}\int_0^\infty\,t^m\,\frac{\mathrm{d}\mu}{t+z}\,,
\qquad m \in \mathbb{N}_0\, ,  \qquad \vert \arg (z) \vert < \pi \, ,  
\end{equation}
in such a way that $F(z)$ takes on the following form:
\begin{equation}
\label{Eq:Stieltjes.7.0}
\displaystyle
F(z)\,=\,\sum_{m=0}^n\,\dfrac{(-)^m}{z^{m+1}}\,\mu_m\,+\,\dfrac{(-)^{n+1}}{z^{n+2}}\,\mu_{n+1}\,\varphi_{n+1}(z)\,.
\end{equation}
From Eq.~(\ref{Eq:Stieltjes.7.0}), it appears that the value of the Stieltjes function $F(z)$ could be retrieved, in principle,
starting from the knowledge of a finite number of terms of the associated Stieltjes series, provided that the corresponding 
converging factor could be estimated in some way.
In particular, if a sufficient approximation to the converging factor could be obtained by a computational approach avoiding the explicit evaluation of the  Stieltjes integral 
into Eq.~(\ref{Eq:Stieltjes.7}), then at least in principle all problems related to the numerical evaluation of the Stieltjes function $F(z)$ via its asymptotic power series (\ref{Eq:Stieltjes.3}) would  definitely be solved, even if this power series converges prohibitively slowly or 
not at all. 

For example, the analysis in \cite{Borghi/Weniger/2015} was ultimately made possible thanks to the following explicit expansion 
obtained for the $n$th-order Euler series converging factor~\cite[Eq.~(52)]{Borghi/2010b},
\begin{equation}
  \label{Eq:BorghiANM}
     {\varphi_{n} ( z )}
       \, = \, \sum_{k=0}^{\infty} \, \frac{L_{k}^{(-1)} (1/z)}{z}
        \, \frac{k!}{(n)_{k+1}} \, ,
\end{equation}
where $L_{k}^{(-1)} (\cdot)$ denotes a generalized Laguerre polynomial
~\cite*[Table 18.3.1]{Olver/Lozier/Boisvert/Clark/2010} and where the symbol $(n)_{k+1}=n(n+1)\ldots(n+k)$ denotes the so-called Pochhammer symbol.
The expansion in Eq. \eqref{Eq:BorghiANM} is just an example of \emph{inverse factorial series} in the discrete index
$n \in \mathbb{N}_{0}$, 
which greatly facilitates the application of finite difference operators, on which Levin-type nonlinear transformations are 
ultimately based. For reader's convenience, the basic definitions and properties of factorial series are briefly recalled in Appendix \ref{App:FactorialSeries}.
A more extensive review can be found, as previously said, in \cite[Appendix B]{Borghi/Weniger/2015}.

It is a remarkably simple task to prove that  the converging factors $\varphi_m(z)$ 
must satisfy a first-order linear recurrence rule. To this end, it is sufficient to start from the integral representation of $\varphi_{m+1}$ given
  in Eq.~(\ref{Eq:Stieltjes.7}), i.e.,  
\begin{equation}
\label{Eq:Stieltjes.7.1}
\displaystyle
\varphi_{m+1}(z) \,=\,\frac {1}{\mu_{m+1}}\int_0^\infty\,t^{m+1}\,\frac{\mathrm{d}\mu}{t+z}\,,
\end{equation}
which can be recast as follows:
\begin{equation}
\label{Eq:Stieltjes.7.2}
\begin{array}{l}
\displaystyle
\varphi_{m+1}(z) \,=\,\frac {1}{\mu_{m+1}}\,\int_0^\infty\,t^{m}\,\frac{t}{t+z}\,\mathrm{d}\mu 
\,=\,
\displaystyle
\frac {1}{\mu_{m+1}}\int_0^\infty\,t^{m}\,\left(1\,-\,\frac{z}{t+z}\right)\,\mathrm{d}\mu\,,
\end{array}
\end{equation}
so that, after taking Eqs.~(\ref{Eq:Stieltjes.0.1}) and~(\ref{Eq:Stieltjes.7}) into account, we have
\begin{equation}
\label{Eq:Stieltjes.7.3}
\displaystyle
\varphi_{m+1} \,=\,
\frac{\mu_m}{\mu_{m+1}}\,\left(1\,-\,z\,\varphi_m\right)\,,\qquad m \ge 0\,,
\end{equation}
Q.E.D.

Accordingly, the converging factor sequence $\{\varphi_m\}_{m=0}^\infty$ must
depend only on $z$ and on the moment ratio sequence $\left\{\dfrac{\mu_m}{\mu_{m+1}}\right\}_{m=0}^\infty$,
a fact which could hardly be obtained in terms of the sole integral representation into Eq. \eqref{Eq:Stieltjes.7.1}.
In particular, the limiting value $\varphi_\infty=\displaystyle\lim_{m\to\infty}\,\varphi_m$
can be extracted directly from Eq.~(\ref{Eq:Stieltjes.7.3}),
 which yields
\begin{equation}
\label{Eq:FactorialExpansionConvergingFactor.2}
\displaystyle
\varphi_{\infty} \,=\,
\frac {\rho_0}{1+\,z\,\rho_0}\,,
\end{equation}
where 
\begin{equation}
\label{Eq:FactorialExpansionConvergingFactor.3}
\displaystyle
\rho_0\,=\,\lim_{m\to\infty}\,\frac{\mu_m}{\mu_{m+1}}\,,
\end{equation}
is directly related to the series convergence radius.

\section{Inverse Factorial series and Weniger transformation}
\label{SubSec:FactorialSeriesWenigerTransformation}

That Stieltjes series converging factors have to satisfy a recurrence relationship, could represent the key 
to understand why the Weniger transformation turns out to be particularly fit for summing series of Stieltjes type. 

Given the $n$th-order converging factor $\varphi_n(z)$, suppose there exist a sequence, say $\{c_k(z)\}_{k=0}^\infty$, whose terms be {\em independent of} $n$, 
such that the following representation holds:
\begin{equation}
\label{Eq:Stieltjes.7.0.1}
\displaystyle
\varphi_n(z)\,=\,\sum_{k=0}^\infty\,
\dfrac{k!}{(n)_{k+1}}\,c_k(z)\,.
\end{equation}
Then, on substituting from Eq.~(\ref{Eq:Stieltjes.7.0.1}) into Eq.~(\ref{Eq:Stieltjes.7.0}), we would have
\begin{equation}
\label{Eq:Stieltjes.7.0.2}
\displaystyle
F(z)\,=\,\sum_{m=0}^n\,a_m(z)\,+\,
a_{m+1}(z)\,
\sum_{k=0}^\infty\,\dfrac{k!}{(n+1)_{k+1}}\,c_k(z)\,,
\end{equation}
where, for brevity, the quantity $a_m(z)=(-)^m\,\mu_m/z^{m+1}$ denotes the $m$th-order term of the Stieltjes series~(\ref{Eq:Stieltjes.3}). 

Suppose for a moment that the factorial series in Eq.~(\ref{Eq:Stieltjes.7.0.1}) consists only of a \emph{finite} number of terms, say $N$.
Then, it is easy to show that the correct value of $F(z)$ can be retrieved within a \emph{finite} number of steps and, more importantly,
\emph{without explicitly knowing} the sequence $\{c_k(z)\}^\infty_{k=0}$ according to the following steps:
\begin{enumerate}

\item Equation~(\ref{Eq:Stieltjes.7.0.2}) is first recast in the form
\begin{equation}
\label{Eq:Stieltjes.7.0.2.1}
\displaystyle
\sum_{k=0}^N\,\dfrac{k!}{(n+1)_{k+1}}\,c_k(z)\,=\,
\dfrac 1{a_{n+1}(z)}\,F(z)\,-\,\dfrac{s_n(z)}{a_{n+1}(z)}\,,
\end{equation}
where $s_n(z)=\displaystyle\sum_{m=0}^n\,a_m(z)$ denotes the $n$th-order partial sum of the series. 

\item Multiply both sides of Eq.~(\ref{Eq:Stieltjes.7.0.2.1}) by $(n+1)_{N+1}$, which gives
\begin{equation}
\label{Eq:Stieltjes.7.0.2.2}
\displaystyle
(n+1)_{N+1}\,\sum_{k=0}^N\,\dfrac{k!}{(n+1)_{k+1}}\,c_k(z)\,=\,
\dfrac {(n+1)_{N+1}}{a_{n+1}(z)}\,F(z)\,-\,\dfrac{(n+1)_{N+1}\,s_n(z)}{a_{n+1}(z)}\,.
\end{equation}

\item Now, the left side of Eq.~(\ref{Eq:Stieltjes.7.0.2.2}) is, without doubts, a $N$th-order polynomial with respect to the variable
$n$. This, in turn, implies that it can be annihilated simply by applying  $N+1$ times the 
first forward difference operator with respect to $n$, defined as 
\begin{equation}
\label{Eq:Stieltjes.7.0.2.3}
\displaystyle
\Delta f_n\,=\,f_{n+1}\,-\,f_{n}\,.
\end{equation}
Accordingly, we have
\begin{equation}
\label{Eq:Stieltjes.7.0.2.5}
\displaystyle
\Delta^{N+1}\,\left\{ (n+1)_{N+1}\,\sum_{k=0}^N\,\dfrac{k!}{(n+1)_{k+1}}\,c_k(z)\right\}\,=\,0\,,
\end{equation}
where the iterated $\Delta^k$ operator turns out to be
\begin{equation}
\label{Eq:Stieltjes.7.0.2.4}
\displaystyle
\Delta^k f_n\,=\,(-1)^k\,\sum_{j=0}^k\,(-)^j\,{{k}\choose{j}}\, f_{n+j}\,.
\end{equation}

\item Finally, substitution from Eq.~(\ref{Eq:Stieltjes.7.0.2.5}) into Eq.~(\ref{Eq:Stieltjes.7.0.2.2}) eventually gives 
\begin{equation}
\label{Eq:Stieltjes.7.0.2.6}
\displaystyle
F(z)\,=\,
\dfrac
{
\Delta^{N+1}\,\left
\{\dfrac{ (n+1)_{N+1}\,s_n(z)}{a_{n+1}(z)}
\right\}
}
{\Delta^{N+1}\,\left\{\dfrac{ (n+1)_{N+1}}{a_{n+1}(z)}\right\}}\,,
\end{equation}
which proves the initial thesis.
\end{enumerate}
It must be remarked again that Eq.~(\ref{Eq:Stieltjes.7.0.2.6}) represents an \emph{exact} result only if the factorial expansion
into Eq.~(\ref{Eq:Stieltjes.7.0.1}) reduces to a \emph{finite} inverse factorial sum. If not, the quantity defined in the 
right side of Eq.~(\ref{Eq:Stieltjes.7.0.2.6}) could be interpreted as the $N$th-order  term of an infinite sequence, say $\{\delta^{(n)}_m\}^\infty_{m=1}$,
which is built up from the partial sum sequence, according to
\begin{equation}
\label{Eq:Stieltjes.7.0.2.6.1}
\displaystyle
\delta^{(n)}_m\,=\,
\dfrac
{
\Delta^{m+1}\,\left
\{\dfrac{ (n+1)_{m+1}\,s_n}{\Delta s_n}
\right\}
}
{\Delta^{m+1}\,\left\{\dfrac{ (n+1)_{m+1}}{\Delta s_n}\right\}}\,,\qquad\qquad m\,\ge\,1\,,n\,\ge\,0\,,
\end{equation}
where $\Delta s_n\,=\,a_{n+1}$ and where, for simplicity, any dependence on the varianble $z$ has to tacitly be intended.
 Equation~(\ref{Eq:Stieltjes.7.0.2.6.1}) defines the so-called  \emph{Weniger transformation} of the partial sum sequence $\{s_m\}^\infty_{m=0}$ \cite{Weniger/1989}. In this case, what would be hopefully desiderable is that, as $m\to\infty$, the $\delta$ sequence defined in
Eq.~(\ref{Eq:Stieltjes.7.0.2.6.1}) could converge, in some sense, to the correct value $F(z)$, i.e.,
\begin{equation}
\label{Eq:Stieltjes.7.0.2.6.1.1}
\displaystyle
\lim_{m\to\infty}\,\delta^{(n)}_m\,=\,F(z)\,, 
\end{equation}
regardless the value of the parameter $n$.

Presently, it does not exist a general convergence theory who could validate, for a typical Stieltjes series, Eq. \eqref{Eq:Stieltjes.7.0.2.6.1.1}.
However, it seems reasonable that the existence of such a theory, assuming it exists, should be based on the representation of the convergence factor given in \eqref{Eq:Stieltjes.7.0.1}. In the next section, an algorithm aimed at building up such a representation will be conceived.

\section{On the Factorial Series Expansion of Stieltjes Series Converging Factor}
\label{SubSec:FactorialSeriesExpansion}

Our idea consists in letting the solution of the difference equation~(\ref{Eq:Stieltjes.7.3}) in the form 
of an inverse factorial series, as follows:
\begin{equation}
\label{Eq:FactorialExpansionConvergingFactor.1}
\displaystyle
\varphi_{m}(z) \,=\,
\varphi_\infty\,+\,\sum_{k=0}^\infty\,
\frac {k!}{(m+\beta)_{k+1}}\,c_k(z)\,,
\end{equation}
where $\varphi_\infty$ is given in Eq.~(\ref{Eq:FactorialExpansionConvergingFactor.2}), $\beta$ is a real nonnegative parameter, and 
$\{c_k(z)\}^\infty_{k=0}$ denotes a sequence {\em independent of} $m$.\footnote{Clearly, it is expected the sequence $\{c_k\}^\infty_{k=0}$ to depend also on $\beta$.
However, for simplicity such a dependence has not been made explicit.} For the moment, it will be assumed $\beta=0$.

Suppose now that also the sequence $\{{\mu_m}/{\mu_{m+1}}\}^\infty_{m=0}$ admits itself an inverse factorial expansion of the form
\begin{equation}
\label{Eq:FactorialExpansionConvergingFactor.1.1}
\displaystyle
\dfrac{\mu_m}{\mu_{m+1}}\,=\,\rho_0\,+\,
\sum_{k=0}^\infty\,
\frac {k!}{(m)_{k+1}}\,b_k\,,
\end{equation}
with $\{b_k\}_{k=0}^\infty$ also being a numerical sequence \emph{independent of} $m$ too.
As we shall see in a moment, Eq.~(\ref{Eq:FactorialExpansionConvergingFactor.1.1}) can be taken as valid under 
rather general conditions.
For example, consider again the Euler series,
 which is the asymptotic series~(\ref{Eq:Stieltjes.3}) with $F(z)=\mathcal{E}(z)$, where
\begin{equation}
\label{Eq:Euler.1}
\displaystyle
\mathcal{E}(z) \, = \, \int_0^\infty\,\frac{\exp(-t)}{t+z}\,\mathrm{d}t\,,
\end{equation}
defines the so-called Euler integral. 
Euler series is a Stieltjes series, with $\mu(t)=1 \, - \, \exp(-t)$
and $\mu_m\,=\,m!$, so that
%
\begin{equation}
\label{Eq:FactorialExpansionConvergingFactor.1.2}
\displaystyle
\dfrac{\mu_m}{\mu_{m+1}}\,=\,
\dfrac 1{m+1}\,=\,
\dfrac 1m\,-\,\dfrac 1m\,+\dfrac 1{m+1}\,=\,
\dfrac 1m\,-\,\dfrac 1{m(m+1)}\,.
\end{equation}
Equation~(\ref{Eq:FactorialExpansionConvergingFactor.1.2})  is of the form~(\ref{Eq:FactorialExpansionConvergingFactor.1.1}) with
$\rho_0=0$, 
$b_0=1$, $b_1=-1$, and $b_k=0$ for $k>1$.  
It is worth noting that the factorial series expansion into Eq.~(\ref{Eq:FactorialExpansionConvergingFactor.1.1})
can always be obtained, in principle, on assuming the moment ratio sequence to admit an inverse  power series representation 
with respect to the index $m$, 
\begin{equation}
\label{Eq:Generalization.1}
\displaystyle
\dfrac{\mu_m}{\mu_{m+1}}\,=\,\sum_{k=0}^\infty\,\dfrac{\rho_k}{m^k}\,=\,\rho_0\,+\,\sum_{k=1}^\infty\,\dfrac{\rho_{k}}{m^k}\,=\,
\rho_0\,+\,\sum_{k=0}^\infty\,\dfrac{\rho_{k+1}}{m^{k+1}}\,,\qquad m \ge 1\,.
\end{equation}
Last series expansion can always be transformed into an inverse factorial series, for example by using the following conversion
formula~\cite{Weniger/2010b}:
\begin{equation}
\label{Eq:Generalization.4}
\displaystyle
\dfrac{\mu_m}{\mu_{m+1}}\,=\,\rho_0\,+\,
\sum_{k=0}^\infty\,\dfrac{\rho_{k+1}}{m^{k+1}}\,=\,
\rho_0\,+\,
\sum_{k=0}^\infty\,\dfrac{k!}{(m)_{k+1}}\,b_k\,,
\end{equation}
where
\begin{equation}
\label{Eq:Generalization.5}
\displaystyle
b_k\,=\,\dfrac{(-1)^k}{k!}\,\sum_{\mu=0}^k\,(-1)^\mu\,
{s}(k,\mu)\,\rho_{\mu+1}\,,
\end{equation}
and the symbol $s(k,\mu)$ denotes the Stirling number of the first order.\footnote{Mathematica Notation for the Stirling number: $s\to \mathrm{S1}$.}
In order to solve Eq.~(\ref{Eq:Stieltjes.7.3}), the following formula will also play a key role:
\begin{equation}
\label{Eq:FactorialExpansionConvergingFactor.1.4}
\begin{array}{l}
\displaystyle
\varphi_{m+1} \,=\,
\varphi_\infty\,+\,
\dfrac{c_0}m\,+\,
\sum_{k=1}^\infty\,
\frac {k!}{(m)_{k+1}}\,(c_k\,-\,c_{k-1})\,,
\end{array}
\end{equation}
whose derivation is detailed in Appendix~\ref{Sec:AppA}.
Another formula which will reveal to be of pivotal importance is that concerning the {\em product} of two factorial series.
In particular, it is not difficult to prove that
\begin{equation}
\label{Eq:FactorialExpansionConvergingFactor.1.5}
\begin{array}{l}
\displaystyle
\left(
\sum_{k=0}^\infty\,
\frac {k!}{(m)_{k+1}}\,b_k
\right)\,
\left(
\sum_{k=0}^\infty\,
\frac {k!}{(m)_{k+1}}\,c_k
\right)\,=\,
\sum_{k=1}^\infty\,
\frac {k!}{(m)_{k+1}}\,d_k\,,
\end{array}
\end{equation}
with the sequence $\{d_k\}_{k=1}^\infty$ being  defined as follows:
\begin{equation}
\label{Eq:FactorialExpansionConvergingFactor.1.6}
\begin{array}{l}
\displaystyle
d_k\,=\,
\dfrac 1{k!}\,
\sum_{\nu=0}^{k-1}\,
\sum_{\mu=0}^{\nu}\,
\nu!\,(\mu+1)_{k-1-\nu}\,
b_{k-1-\nu}\,c_{\nu-\mu}\,,
\end{array}
\end{equation}
or, equivalently,
\begin{equation}
\label{Eq:FactorialExpansionConvergingFactor.1.7}
\begin{array}{l}
\displaystyle
d_k\,=\,
\dfrac 1{k!}\,
\sum_{\nu=0}^{k-1}\,
\sum_{\lambda=0}^{\nu}\,
\nu!\,(\nu-\lambda+1)_{k-1-\nu}\,
b_{k-1-\nu}\,c_{\lambda}\,.
\end{array}
\end{equation}
Now, all we have to do is to substitute from Eqs.~(\ref{Eq:Generalization.4}),~(\ref{Eq:FactorialExpansionConvergingFactor.1}), 
and~(\ref{Eq:FactorialExpansionConvergingFactor.1.4}) into Eq.~(\ref{Eq:Stieltjes.7.3}). 
Nontrivial mathematical steps, detailed in Appendix~\ref{Sec:AppB}, lead to
\begin{equation}
\label{Eq:FactorialExpansionConvergingFactor.1.9.TER}
\left\{
\begin{array}{l}
\displaystyle
{c_{-1}}\,=\,0\,,\\
\\
{c_0}\,=\,\dfrac{b_0}{(1+z\rho_0)^2}\,,\\
\\
c_k\,=\,\dfrac{c_{k-1}}{1+z\rho_0}\,+\,\dfrac{b_k}{(1+z\rho_0)^2}\,-\,\dfrac{z}{1+z\rho_0}\,d_k\,,\,\qquad k\ge 1\,,
\end{array}
\right.
\end{equation}
which, together with 
Eqs.~(\ref{Eq:FactorialExpansionConvergingFactor.1.1})
and~(\ref{Eq:FactorialExpansionConvergingFactor.1.6}),
represents the main result of the present paper. 

We have proved that the $n$th-order truncation error of a Stieltjes series can be, 
under quite general hypotheses, expressed as the product of the first term of the series neglected by the $n$th-order partial sum,
times an inverse factorial series whose single terms are defined through a triangular recurrence relation involving, apart from $z$, only the asymptotic expansion of the moment ratio sequence $\{\mu_{m}/\mu_{m+1}\}$ for $m\to\infty$.

It should be noted how, in the case of divergent series with $\rho_0=0$, Equation~(\ref{Eq:FactorialExpansionConvergingFactor.1.9.TER}) takes on the following easier form: 
\begin{equation}
\label{Eq:FactorialExpansionConvergingFactor.1.9.QUATER}
\left\{
\begin{array}{l}
\displaystyle
{c_{-1}}\,=\,0\,,\\
\\
{c_0}\,=\,{b_0}\,,\\
\\
c_k\,=\,{c_{k-1}}\,+\,{b_k}\,-\,{z}\,d_k\,,\,\qquad k>1\,.
\end{array}
\right.
\end{equation}

\section{A Few Examples of Applications}
\label{Sec:Examples}
%

\subsection{Revisiting Euler series}
\label{SubSec:ES}

In Ref.~\cite{Borghi/2010b}, the factorial representation~(\ref{Eq:FactorialExpansionConvergingFactor.1}) of the converging factor of the Euler series was  found  for the particular case of $\beta=0$, being
\begin{equation}
\label{Eq:Euler.0.2.1.1}
\displaystyle
\varphi_m(z)\,=\,\sum_{k=0}^\infty\,
\dfrac{k!}{(m)_{k+1}}\,zL^{(-1)}_k(z)\,,
\end{equation}
where the symbol  $L^{(\alpha)}_k(\cdot)$ denotes the associate 
Laguerre polynomial~\cite{Abramowitz/Stegun/1972}. 
Equation~(\ref{Eq:Euler.0.2.1.1}) was obtained by using 
the general approach proposed in~\cite{Weniger/2007a}. 
This led to a nontrivial derivation  involving  the use of Bernoulli as well as exponential
polynomials~\cite{Borghi/2010b}.
In the present section, Eq.~(\ref{Eq:Euler.0.2.1.1}) will now be re-derived by using the results previously obtained
and generalized to arbitrary positive values of $\beta$.  To this end, the first, 
preliminary step is to find the factorial expansion of the moment ratio  in terms of $m+\beta$. 
This can be achieved, for instance, by recasting the moment ratio as follows: 
\begin{equation}
\label{Eq:Euler.1.2.6.1}
\displaystyle
\dfrac{\mu_m}{\mu_{m+1}}\,=\,\dfrac 1{m+1}\,=\,\dfrac 1{m+\beta-(\beta-1)}\,,
\end{equation}
and then by using the fundamental  expansion \cite[Eq.\ (3) on p.\ 77]{Nielsen/1965}:\footnote{In the books by \citet[p.\ 291]{Milne-Thomson/1981} and by \citet[p.\
177]{Paris/Kaminski/2001}, this expansion is called Waring's formula, but
in Nielsen's book \citep[Footnote 2 on p.\ 77]{Nielsen/1965} this formula
is attributed to \citet{Stirling/1730}, and in the book by \citet[p.\
199]{Noerlund/1954} it is attributed to Nicole (compare also \citet[p.\
30]{Tweedie/1917}).
}
\begin{equation}
  \label{WaringFormula}
\frac{1}{z-w} \; = \; \sum_{n=0}^{\infty} \,
\frac{(w)_n}{(z)_{n+1}} \, , \qquad \mathrm{Re} (z-w) > 0 \, ,
\end{equation}
which gives at once
\begin{equation}
\label{Eq:Euler.1.2.6.2}
\displaystyle
\dfrac{\mu_m}{\mu_{m+1}}\,=\,
\sum_{k=0}^\infty\,\dfrac{k!}{(m+\beta)_{k+1}}\,\dfrac{(\beta-1)_k}{k!}\,.
\end{equation}
On comparing Eq.~(\ref{Eq:Euler.1.2.6.2}) with Eq.~(\ref{Eq:FactorialExpansionConvergingFactor.1.1})
and on taking into account that $\rho_0=0$, we  have
\begin{equation}
\label{Eq:Euler.1.2.6.3}
\displaystyle
b_k\,=\,\dfrac{(\beta-1)_k}{k!}\,,\qquad\qquad k\ge 0\,,\quad\beta\ge 0\,.
\end{equation}

It is worth showing an alternative approach, based on the asymptotic expansion in Eq.~(\ref{Eq:Generalization.1}).
In particular, the expanding coefficients $\{\rho_k\}^\infty_{k=0}$ can be found by writing  Eq.~(\ref{Eq:Euler.1.2.6.1})
as follows:
\begin{equation}
\label{Eq:Euler.1.2.6.3.1}
\displaystyle
\dfrac {\mu_m}{\mu_{m+1}}\,=\,\dfrac 1{m+\beta-(\beta-1)}\,=\,
\dfrac 1{m+\beta}\,\dfrac 1{1\,-\,\dfrac {\beta\,-\,1}{m\,+\,\beta}}\,,
\end{equation}
and by expanding the second factor as a geometric series, so that
\begin{equation}
\label{Eq:Euler.1.2.6.3.2}
\displaystyle
\dfrac {\mu_m}{\mu_{m+1}}\,=\,
\sum_{k=0}^\infty\,\dfrac{(\beta-1)^k}{(m+\beta)^{k+1}}\,.
\end{equation}
On comparing Eq.~(\ref{Eq:Euler.1.2.6.3.2}) with Eq.~(\ref{Eq:Generalization.1}), we have $\rho_0=0$ and $\rho_{k}=(\beta-1)^{k-1}$, for $k\ge 1$, so that
Eq.~(\ref{Eq:Generalization.5}) leads to
\begin{equation}
\label{Eq:Euler.1.2.6.3.3}
\displaystyle
b_k\,=\,\dfrac{(-1)^k}{k!}\,\sum_{\mu=0}^k\,(-1)^\mu\,
{s}(k,\mu)\,(\beta-1)^{\mu}\,=\,\dfrac{(\beta-1)^k}{k!}\,,
\end{equation}
%
which proves again Eq.~(\ref{Eq:Euler.1.2.6.3}).

Once the factorial expansion of the moment ratio has been obtained, it is sufficient to substitute from Eq.~(\ref{Eq:Euler.1.2.6.3}) into
Eq.~(\ref{Eq:FactorialExpansionConvergingFactor.1.9.TER}) to obtain, with a small amount of elementary algebra,
\begin{equation}
\label{Eq:Euler.1.2.6.3.4}
\begin{array}{l}
\displaystyle
c_0\,=\,1\,,\\
\\
c_1\,=\,\beta\,-\,z\,=\,L^{(\beta-1)}_1(z)\,,\\
\\
c_2\,=\,\dfrac{1}{2} \left(\beta ^2+\beta +z^2-2 (\beta +1) z\right)\,=\,L^{(\beta-1)}_2(z)\,,
\end{array}
\end{equation}
and so on, for all $k>2$.
Accordingly, Eq. (\ref{Eq:Euler.0.2.1.1}) becomes
\begin{equation}
\label{Eq:Euler.0.2.1.1BIS}
\displaystyle
\varphi_m(z)\,=\,\sum_{k=0}^\infty\,
\dfrac{k!}{(m+\beta)_{k+1}}\,zL^{(\beta-1)}_k(z)\,,\qquad \beta \ge 0\,.
\end{equation}

\subsection{Error function}
\label{Sec:ErrorFunction}

The following  Stieltjes series:
\begin{equation}
\label{Sec:DiscussionOpenProblems.1}
\begin{array}{l}
\displaystyle
\sum_{m=0}^\infty\,\frac{(-1)^m}{z^{m+1}}\,\Gamma\left(m+\frac 12\right)\,,
\end{array}
\end{equation}
where $\Gamma(\cdot)$ denotes the Gamma function~\cite{NIST:DLMF}, is asymptotic to the complementary error function erfc$(\cdot)$~\cite{NIST:DLMF}, being
\begin{equation}
\label{Sec:DiscussionOpenProblems.1.1}
\begin{array}{l}
\displaystyle
\sum_{m=0}^\infty\,\frac{(-1)^m}{z^{m+1}}\,\Gamma\left(m+\frac 12\right)\,\sim\,
\frac{\pi\,\exp(z)\,\mathrm{erfc}(\sqrt z)}{\sqrt z}\,,\qquad z\to\infty\,.
\end{array}
\end{equation}
In order to find the factorial expansion of the corresponding converging factor,
it is sufficient to recast the moment ratio as follows:
\begin{equation}
\label{Sec:DiscussionOpenProblems.2.1}
\begin{array}{l}
\displaystyle
\dfrac{\mu_m}{\mu_{m+1}}\,=\,\dfrac{\Gamma\left(m\,+\,\dfrac 12\right)}{\Gamma\left(m\,+\,\dfrac 32\right)}\,=\,
\dfrac 1{m\,+\,\dfrac 12}\,=\,\dfrac 1{\left(m\,-\,\dfrac 12\right)\,+\,1}\,,
\end{array}
\end{equation}
and to note that is formally identical to Eq.~(\ref{Eq:Euler.1.2.6.1}), once the transformation $m\to m-1/2$ be applied.
Accordingly, it is not difficult to prove that all steps leading to Eq.~(\ref{Eq:Euler.0.2.1.1BIS})
for the Euler series can be repeated for the asymptotic series~(\ref{Sec:DiscussionOpenProblems.1}), which yields
\begin{equation}
\label{Sec:DiscussionOpenProblems.3}
\begin{array}{l}
\displaystyle
\varphi_m(z)\,=\,
\sum_{k=0}^\infty\,
\frac{k!}{\left(m\,+\,\beta\,-\,\frac 12\right)_{k+1}}\,zL^{(\beta-1)}_k(z)\,,\qquad\,m>1\,.
\end{array}
\end{equation}
%


\subsection{The logarithm function}
\label{SubSec:Logarithm}

Another important Stieltjes function is the natural logarithm. 
Consider the function $\mathcal{L}(z)$ defined as follows:
\begin{equation}
\label{Eq:Madelung.1}
\displaystyle
\mathcal{L}(z)\,=\,
\int_0^1\,\frac{\mathrm{d}t}{z\,+\,t}\,=\,\log\left(1\,+\,\frac 1z\right)\,.
\end{equation}
The Stieltjes series of the form~(\ref{Eq:Stieltjes.3}) asymptotic to $\mathcal{L}(z)$, is characterized by the following moment sequence:
\begin{equation}
\label{Eq:Madelung.2}
\displaystyle
\mu_m\,=\,\frac 1{m+1}\,,\qquad\,m\ge 0\,,
\end{equation}
so that 
\begin{equation}
\label{Eq:Madelung.2.1}
\displaystyle
\dfrac{\mu_m}{\mu_{m+1}}\,=\,\frac {m+2}{m+1}\,=\,1\,+\,\dfrac 1{m+1}\,.
\end{equation}
Similarly as happened for the error function, also the factorial expansion of the converging factor of the logarithm
can be obtained in a surprisingly simple way. In fact, it is sufficient to compare Eq.~(\ref{Eq:Madelung.2.1}) with Eq.~(\ref{Eq:Generalization.4}), written with $m+\beta$ in place of $m$, and on taking Eq.~(\ref{Eq:Euler.1.2.6.2}) into account, to have
\begin{equation}
\label{Eq:Madelung.2.2}
\displaystyle
\dfrac {\mu_m}{\mu_{m+1}}\,=\,1\,+\,
\sum_{k=0}^\infty\,\dfrac{k!}{(m+\beta)^{k+1}}\,\dfrac{(\beta-1)^k}{k!}\,,
\end{equation}
from which it follows that
\begin{equation}
\label{Eq:Madelung.2.3}
\begin{array}{l}
\displaystyle
\rho_0\,=\,1\,\\
\\
b_k\,=\,\dfrac{(\beta-1)^k}{k!}\,.
\end{array}
\end{equation}
Then, on substituting from Eq.~(\ref{Eq:Madelung.2.3}) into Eq.~(\ref{Eq:FactorialExpansionConvergingFactor.1.9.TER})
after simple algebra we obtain
\begin{equation}
\label{Eq:Madelung.2.4}
\begin{array}{l}
\displaystyle
c_0\,=\,\frac 1{(1+z)^2}\,,\\
\\
\displaystyle
c_1\,=\,\frac {\beta\,+\, (\beta-2)\,z}{(1+z)^3}\,,\\
\\
\displaystyle
c_2\,=\,\frac { \beta (1+ \beta)\,+\, 2 (\beta-2) (\beta+1)\,z\,+\, (\beta-2) (\beta-1)\,z^2}{2\,(1+z)^4}\,,
\end{array}
\end{equation}
and so on, for all $k>2$.
Remarkably, it is also possible to give Eq.~(\ref{Eq:Madelung.2.4}) the following closed-form expression:
%
\begin{equation}
\label{Eq:Madelung.2.3.2}
\displaystyle
c_k(z)\,=\,\frac{(-1)^k}{(1+z)^2}\,
P_k^{(2-\beta-k,\beta-1)}\left(\frac{z-1}{z+1}\right)\,,\qquad\,k\ge 0\,,
\end{equation}
where the symbol $P_n^{(\alpha,\beta)}(\cdot)$ denotes the Jacobi polynomial~\cite{NIST:DLMF}.

\subsection{Lerch's trascendental function}
\label{Sec:Lerch}

The logarithm function is a particular case of a more general Stieltjes function, namely the Lerch trascendental~\cite{NIST:DLMF},
whose asymptotic Stieltjes series~(\ref{Eq:Stieltjes.3}) is characterized by the following moment sequence:
\begin{equation}
\label{Eq:Lerch.1}
\displaystyle
\mu_m\,=\,\frac 1{(m+\alpha)^s}\,,\qquad\,m\ge 0\,,
\end{equation}
where, for simplicity, it will be assumed $\alpha>0$ and $s>1$.
The moment ratio then turns out to be
\begin{equation}
\label{Eq:Lerch.2}
\displaystyle
\dfrac{\mu_m}{\mu_{m+1}}\,=\,
\left(1\,+\,\dfrac 1{m+\alpha}\right)^s\,=\,
\left(1\,+\,\dfrac 1{m+\beta\,+\,(\alpha-\beta)}\right)^s\,,\qquad\qquad m \ge 0\,,
\end{equation}
whose inverse factorial expansion can be obtained by using Eqs.~(\ref{Eq:Generalization.1})\,-\,(\ref{Eq:Generalization.5}).
To this end, consider first the case $\beta=0$.  Then, the asymptotic inverse power expansion with respect to $m$ gives
\begin{equation}
\label{Eq:Lerch.3}
\begin{array}{l}
\displaystyle
\rho_0\,=\,1\,,\\
\\
\rho_1\,=\,s\,,\\
\\
\displaystyle
\rho_2\,=\,\frac 12\,s(s-1) \,-\,  s \alpha \,,\\
\\
\displaystyle
\rho_3\,=\,\frac 16\,s(s-1) (s-2)  - s(s-1) \alpha + s \alpha^2\,,\\
\\
\displaystyle
\rho_4\,=\,\frac 1{24}\,s(s-1) (s-2) (s-3)  \,-\,\frac 12\,s(s-1) (s-2) \alpha \,+\,
\frac 32 \,s(s-1) \alpha^2\,-\,s \alpha^3\,,\\
\\
\ldots
\end{array}
\end{equation}
Accordingly, on using again Eq.~(\ref{Eq:Generalization.5}), the expanding coefficients $b_k$ of the factorial expansion
of  $\dfrac{\mu_m}{\mu_{m+1}}$ turn out to be
\begin{equation}
\label{Eq:Lerch.4}
\begin{array}{l}
\displaystyle
b_0\,=\,s\,,\\
\\
\displaystyle
b_1\,=\,\frac 12\,s (s-1 - 2 \alpha) \,,\\
\\
\displaystyle
b_2\,=\,\frac 1{12}\,s (s^2-1 - 6 s \alpha + 6  \alpha^2) \,,\\
\\
\displaystyle
b_3\,=\,\frac 1{144}\,s\,( s (s (6 + s)-1) - 12 s (s+3) \alpha \,+\,  36 ( s+1) \alpha^2 - 24 \alpha^3\,-\,6)\,,\\
\\
\ldots
\end{array}
\end{equation}
and our recursive algorithm provides the following analytical expressions of the first terms of the sequence $\{c_k\}^\infty_{k=0}$:
\begin{equation}
\label{Eq:Lerch.5}
\begin{array}{l}
\displaystyle
c_0\,=\,\dfrac s{(1 + z)^2}\,,\\
\\
c_1\,=\,-\,\dfrac{(s ((s+1) (z-1) + 2 (1+z) \alpha)}{2 (1 + z)^3}\,,\\
\\
c_2\,=\,\dfrac{s ((1 + s) (5 + s + s ( z-4) z - z (8 + z)) + 
   6 ( s ( z-1)-2) (1 + z) \alpha + 6 (1 + z)^2 \alpha^2)}{12 (1 + z)^4}\,,\\
\\
\ldots
\end{array}
\end{equation}
Note that, in order to recover the more general factorial expansion of the converging factor~(\ref{Eq:FactorialExpansionConvergingFactor.1}) valid
for a typical $\beta>0$, it is sufficient,  according to Eq.~(\ref{Eq:Lerch.2}), to replace $\alpha$ by $\alpha-\beta$ into Eq.~(\ref{Eq:Lerch.5}). 
Differently from Eq. (\ref{Eq:Madelung.2.4}), we were not able to give Eq. (\ref{Eq:Lerch.5}) a closed form expression.

\section{Conclusions}
\label{Sec:Conclusions}

Although the practical usefulness of divergent series has always been more important than concerns about mathematical rigor, 
especially in physics, a theoretical understanding of the mathematical machinery behind Levin-type sequence transformations should be 
highly desiderable. An important achievement in this direction has been established in~\cite{Borghi/Weniger/2015}, 
where it was rigorously proved that a particular type of Levin-type transformation, the Weniger transformation~\cite{Weniger/1989},
is not only able to resum the Euler series,~but also that it turns out to be ``exponentially faster''  than Pad\'e approximants
in accomplishing such a task. Euler series is the paradigm of factorial divergence. Several other series expansions occurring in mathematical
physics are 
factorially or even hyperfactorially diverging. The unquestionable superiority of the Weniger transformation over Pad\'e in resumming the Euler series proved in~\cite{Borghi/Weniger/2015} were not limited to such a specific case, as it was shown in the past, although only on a purely numerical basis.
Levin-type transformations in general, and Weniger's transformation in particular, are able to dramatically accelerate
the convergence of slowly convergent series as well as to resum wildly divergent series, some of them even not resummable 
via Pad\'e. The retrieving capabilities of Weniger's transformation are strictly related to the possibility of representing the 
series converging factor of any finite order as a suitable {inverse factorial series}.

In the present paper, which is the natural continuation of Ref. \cite{Borghi/Weniger/2015}, it has been proved that the 
$n$th-order converging factor of a Stieltjes series can be expressed through an inverse factorial series, whose terms are uniquely determined, starting from the knowledge of  the series moment ratio, through a linear recurrence relation. Our constructive proof has the form of 
an algorithm, developed starting from the first-order difference equation which has to be satisfied by the convergent factors of {\em any} Stieltjes series. 
A few examples of applications of our algorithm have also presented, in order to show its effectiveness and implemementation ease.

The class of Stieltjes series is of fundamental interest in mathematical physics. 
Recently, it was rigorously shown that the factorially divergent perturbation expansion for the energy eigenvalue of the $\mathcal{P T}$-symmetric Hamiltonian $H(\lambda) = p^2 + \frac{1}{4} x^2 + \mathrm{i} \lambda x^3$ is a Stieltjes series~\cite{Grecchi/Maioli/Martinez/2009,Grecchi/Martinez/2013}, thus proving a conjecture formulated  ten years before by Bender and Weniger~\cite{Bender/Weniger/2001}. 
A much more recent result is the proof that also the character of the celebrated Bessel solution of Kepler's equation \cite{Colwell/1993} is Stieltjes \cite{Borghi/2024}. 
Weniger's transformation turned out to be successful, in the past, for resumming extremely violently divergent perturbation expansions
\cite{Weniger/1996c,Weniger/1996e,%
  Weniger/Cizek/Vinette/1991,Weniger/Cizek/Vinette/1993}, where 
Pad\'{e} approximants were not powerful enough to achieve anything substantial,
as in the case of Rayleigh-Schr\"{o}inger perturbation series for the sextic anharmonic oscillator,
or even not able to sum, as for the more challenging octic case \cite{Graffi/Grecchi/1978}. 
There is no doubt that a better theoretical understanding of the summation of
hyper-factorially divergent expansions via Weniger's transformation would also be highly desirable.

EJW and I were convinced that the results obtained in \cite{Borghi/Weniger/2015} could have been the first, preliminary step toward the development of a general convergence theory of the summation of Stieltjes series with Levin-type transformations in place of Pad\'e approximants. What has been shown in the present paper might be, in my opinion, a small leap in this direction. 

\appendix

\section[\appendixname~\thesection]{Basic Properties of Factorial Series}
\label{App:FactorialSeries}

Let $\Omega (z)$ be a function that vanishes as $z \to + \infty$. A
\emph{factorial series} for $\Omega (z)$ is an expansion in terms of
inverse Pochhammer symbols of the following kind:
\begin{equation}
  \label{DefFactSer}
  \Omega (z) \; = \; 
   \frac {b_{0}} {z} + \frac {b_{1} 1!} {z(z+1)}
    + \frac {b_{2} 2!} {z(z+1)(z+2)} + \cdots \; = \; 
     \sum_{\nu=0}^{\infty} \frac {b_{\nu} {\nu}!} {(z)_{\nu+1}} \, .
\end{equation}
In general, $\Omega (z)$ will have simple poles at $z = - m$ with $m \in
\mathbb{N}_0$.

Factorial series had been employed already in Stirling's classic book
\emph{Methodus Differentialis} \citep{Stirling/1730,Stirling/1749}, which
appeared in \citeyear{Stirling/1730}. 

In the nineteenth and the early twentieth century, the theory of
factorial series was fully developed by a variety of authors. Fairly
complete surveys of the older literature as well as thorough treatments
of their properties can be found in books on finite differences by
\citet{Meschkowski/1959}, \citet{Milne-Thomson/1981},
\citet{Nielsen/1965}, and
\citet{Noerlund/1926,Noerlund/1929,Noerlund/1954}. Factorial series are
also discussed in the books by \citet{Knopp/1964} and
\citet{Nielsen/1909} on infinite series. 

Compared to inverse power series, factorial series possess vastly
superior convergence properties . If we use (see for example \citet*[Eq.\
(5.11.12 )]{Olver/Lozier/Boisvert/Clark/2010})
\begin{equation}
  \label{AsyGammaRatio}
\Gamma(z+a)/\Gamma(z+b) \; = \; z^{a-b} \,
\bigl[ 1 + \mathrm{O} (1/z) \bigr] \, , \qquad z \to \infty \, ,
\end{equation}
we obtain the asymptotic estimate $n!/(z)_{n+1} = \mathrm{O} (n^{-z})$ as
$n \to \infty$. Thus, the factorial series (\ref{DefFactSer}) converges
with the possible exception of the points $z = - m$ with
$m \in \mathbb{N}_0$ if and only if the associated Dirichlet series
$\bar{\Omega} (z) = \sum_{n=1}^{\infty} b_{n}/n^z$ converges (see for
example \cite[p.\ 262]{Knopp/1964} or \cite[p.\ 167]{Landau/1906}).
Accordingly, a factorial series converges for sufficiently large
$\Re (z)$ even if the reduced series coefficients $b_{n}$ in
(\ref{DefFactSer}) grow like a fixed power $n^{\alpha}$ with
$\Re (\alpha) > 0$ as $n \to \infty$. Thus, a factorial series can
converge even if its effective coefficients $b_{n} n!$ diverge quite
wildly as $n \to \infty$.

Factorial series can also be expressed in terms of the beta function
$B (x, y) = \Gamma(x) \Gamma(y) / \Gamma(x+y)$ \citep*[Eq.\
(5.12.1)]{Olver/Lozier/Boisvert/Clark/2010}. Because of
$B (z, n+1) = n!/(z)_{n+1}$, the factorial series in (\ref{DefFactSer})
can be interpreted as an expansion in terms of special beta functions
(see for instance \citet[p.\ 288]{Milne-Thomson/1981} or \citet[p.\
175]{Paris/Kaminski/2001}):
\begin{equation}
  \label{FactSerBetaExpan}
  \Omega (z) \; = \; \sum_{n=0}^{\infty} \, b_{n} \, B (z, n+1) \, .
\end{equation}
The beta function possesses the integral representation (see for example
\citep*[Eq.\ (5.12.1)]{Olver/Lozier/Boisvert/Clark/2010})
\begin{equation}
  B (x, y) \; = \; \int_{0}^{1} \, t^{x-1} \, (1-t)^{y-1} \, \mathrm{d} t 
   \, , \qquad \Re (x), \Re (y) > 0 \, ,
\end{equation}
which implies
\begin{equation}
  \label{IntRepFSterm}
  B (z, n+1) \; = \; \frac{n!}{(z)_{n+1}} \; = \;
  \int_{0}^{1} \, t^{z-1} \, (1-t)^{n} \, \mathrm{d} t
  \, , \qquad \Re (z) > 0 \, \quad n \in \mathbb{N}_0 \, .
\end{equation}
If we insert this into (\ref{FactSerBetaExpan}) and interchange
integration and summation, we obtain the following very important
integral representation \citep[Sec. I on p.\ 244]{Nielsen/1965}:
\begin{subequations}
  \label{FS_IntRep}
  \begin{align}
    \label{FS_IntRep_a}
    \Omega (z) & \; = \; \int_{0}^{1} \, t^{z-1} \, \varphi_{\Omega} (t)
     \, \mathrm{d} t \, , \qquad \Re (z) > 0 \, ,
    \\
    \label{FS_IntRep_b}
    \varphi_{\Omega} (t) & \; = \; 
     \sum_{n=0}^{\infty} \, b_{n} \, (1-t)^{n} \, .
  \end{align}
\end{subequations}
Often, the properties of $\Omega (z)$ can be studied more easily via this
integral representation than via the defining factorial series
\eqref{DefFactSer} (see for example \citep[Ch. XVII]{Nielsen/1965}).

For our purposes, the mathematical definition \eqref{DefFactSer} is not
completely satisfactory and it is more natural to use instead a slightly
different definition. Let $\tilde{\Omega} (z)$ be a function which is
finite, but non-zero as $z \to +\infty$, and let us assume that it
possesses a slightly more general factorial series expansion of the
following kind:
\begin{equation}
  \label{DefFactSer_mod}
  \tilde{\Omega} (z) \; = \; \beta_{0} + 
   \frac {\beta_{1}} {z} + \frac {\beta_{2}} {z(z+1)}
    + \frac {\beta_{3}} {z(z+1)(z+2)} + \cdots \; = \; 
     \sum_{\nu=0}^{\infty} \frac {\beta_{\nu}} {(z)_{\nu}} \, .
\end{equation}
Obviously, we have $\lim_{z \to +\infty} \tilde{\Omega} (z) =
\beta_{0}$. We obtain the mathematical definition \eqref{DefFactSer} 
if we set in (\ref{DefFactSer_mod})
 $\beta_{0} = 0$ and
$\beta_{\nu+1} = b_{\nu} {\nu}!$.

\section[\appendixname~\thesection]{Proof of Eq.~(\ref{Eq:FactorialExpansionConvergingFactor.1.4})}
\label{Sec:AppA}

First of all, $\varphi_{m+1}$ is recast as follows:
\begin{equation}
\label{Eq:appA.0}
\displaystyle
\varphi_{m+1}\,=\,\varphi_{m}\,+\,\Delta\,\varphi_{m} \,.
\end{equation}
Then, on applying the forward difference operator $\Delta$ to both sides of Eq.~(\ref{Eq:FactorialExpansionConvergingFactor.1}),
we have
\begin{equation}
\label{Eq:appA.1}
\displaystyle
\Delta\,\varphi_{m} \,=\,
\Delta\,\sum_{k=0}^\infty\,
\frac {k!}{(m)_{k+1}}\,c_k\,=\,
\sum_{k=0}^\infty\,
k!\,c_k\,\Delta \left(\frac {1}{(m)_{k+1}}\right)\,,
\end{equation}
where (see~\cite[Eq.~(13.2-12)]{Weniger/1989})
\begin{equation}
\label{Eq:appA.2}
\displaystyle
\Delta \left(\frac {1}{(m)_{k+1}}\right)\,=\,
-\frac {k+1}{(m)_{k+2}}\,,
\end{equation}
so that
\begin{equation}
\label{Eq:appA.3}
\displaystyle
\Delta\,\varphi_{m} \,=\,-\sum_{k=1}^\infty\,
\frac {k!}{(m)_{k+1}}\,c_{k-1}\,.
\end{equation}
Then, on taking again Eq.~(\ref{Eq:appA.0}) into account, from Eqs.~(\ref{Eq:FactorialExpansionConvergingFactor.1})
and~(\ref{Eq:appA.3}), Eq.~(\ref{Eq:FactorialExpansionConvergingFactor.1.4}) follows after simple rearranging.

\section[\appendixname~\thesection]{Proof of Eq.~(\ref{Eq:FactorialExpansionConvergingFactor.1.9.TER})}
\label{Sec:AppB}
%

We start by recasting  Eq.~(\ref{Eq:Stieltjes.7.3}) as follows:
\begin{equation}
\label{Eq:appB.1}
\displaystyle
\varphi_{m+1} \,=\,
\frac{\mu_m}{\mu_{m+1}}\,-\,z\,\frac{\mu_m}{\mu_{m+1}}\,\varphi_m\,,
\end{equation}
which, on taking Eqs.~(\ref{Eq:FactorialExpansionConvergingFactor.1}),~(\ref{Eq:Generalization.4}),
and~(\ref{Eq:FactorialExpansionConvergingFactor.1.4}) into account, can be recast as follows:
\begin{equation}
\label{Eq:appB.2}
\begin{array}{l}
\displaystyle
\varphi_\infty\,+\,
\sum_{k=0}^\infty\,\frac {k!}{(m)_{k+1}}\,\Delta c_{k-1}\,=\, 
\rho_0\,+\,
\sum_{k=0}^\infty\,\dfrac{k!}{(m)_{k+1}}\,b_k\,\\
\\
\displaystyle
-\,z\,
\left(\rho_0\,+\,
\sum_{k=0}^\infty\,\dfrac{k!}{(m)_{k+1}}\,b_k
\right)
\left(\varphi_\infty\,+\,
\sum_{k=0}^\infty\,\dfrac{k!}{(m)_{k+1}}\,c_k
\right)\,,
\end{array}
\end{equation}
where, for simplicity, it has been assumed formally $c_{-1}\,=\,0$.
Equation~(\ref{Eq:appB.2}) can be rearranged, on taking Eq.~(\ref{Eq:FactorialExpansionConvergingFactor.2}) into account,
as follows:
\begin{equation}
\label{Eq:appB.3}
\begin{array}{l}
\displaystyle
\sum_{k=0}^\infty\,\frac {k!}{(m)_{k+1}}\,\Delta c_{k-1}\,=\,
\frac 1{1\,+\,z\rho_0}\,\sum_{k=0}^\infty\,\dfrac{k!}{(m)_{k+1}}\,b_k\,\\
\\
\displaystyle
-\,z\rho_0\,\sum_{k=0}^\infty\,\dfrac{k!}{(m)_{k+1}}\,c_k\,\,-\,
z\,\left(\sum_{k=0}^\infty\,\dfrac{k!}{(m)_{k+1}}\,b_k\right)\,
\left(\sum_{k=0}^\infty\,\dfrac{k!}{(m)_{k+1}}\,c_k\right)\,,
\end{array}
\end{equation}
Finally, on using 
Eq.~(\ref{Eq:FactorialExpansionConvergingFactor.1.5}) into Eq.~(\ref{Eq:appB.3}),
after simple algebra we  obtain
\begin{equation}
\label{Eq:appB.4}
\left\{
\begin{array}{l}
\displaystyle
{c_0}\,=\,\dfrac{b_0}{1+z\rho_0}\,-\,z\rho_0\,c_0\,,\\
\\
\Delta c_{k-1}\,=\,\dfrac{b_k}{1+z\rho_0}\,-\,z\rho_0\,c_k\,-\,z\,d_k\,,\,\qquad k\ge1\,,
\end{array}
\right.
\end{equation}
from which Eq.~(\ref{Eq:FactorialExpansionConvergingFactor.1.9.TER}) follows.

\acknowledgments{I wish to thank Turi Maria Spinozzi for his useful comments and help.}

\end{document}